\documentclass{amsart}
\usepackage{graphicx}

\newtheorem{thm}{Theorem}[section]

\newtheorem{lem}[thm]{Lemma}

\theoremstyle{definition}

\theoremstyle{remark}

\numberwithin{equation}{section}
\newcommand{\bea}{\begin{eqnarray}}
\newcommand{\eea}{\end{eqnarray}}

\def\s{\sigma}

\def\t{\tau}

\def\M{\mathcal{M}}
\def\F{\mathcal{F}}

\def\m{\mu}

\def\L{\Lambda}

\def\e{\varepsilon}

\def\bb{{\mathbb B}}
\def\bc{{\mathbb C}}

\def\bk{{\mathbb K}}

\def\bn{{\mathbb N}}

\def\gb{{\mathbf{g}}}
\def\zb{{\mathbf{z}}}
\def\kb{{\mathbf{k}}}
\def\nb{{\mathbf{n}}}

\def\pb{{\mathbf{p}}}
\def\Nb{{\mathbf{N}}}
\def\Tb{{\mathbf{T}}}

\begin{document}

\title[  multiparameter WEIGHTED ergodic theorem ] {
on multiparameter WEIGHTED ergodic theorem for Noncommutative
$L_{p}$-spaces}

\author{Farrukh Mukhamedov}
\address{Farrukh Mukhamedov\\
Departamento de Fisica,
Universidade de Aveiro\\
Campus Universitario de Santiago\\
3810-193 Aveiro, Portugal} \email{{\tt far75m@yandex.ru}, {\tt
farruh@fis.ua.pt}}

\author{Maksut Mukhamedov}
\address{Maksut  Mukhemedov\\
Faculty of Economy \\
Department of High Mathematics\\
Tashkent State Agrar University\\
Tashkent, Uzbekistan}
\author{Seyit Temir}
\address{Seyit Temir\\
Department of Mathematics\\
Art and Science Faculty\\
Harran University, 63200, Sanliurfa, Turkey} \email{{\tt
temirseyit@harran.edu.tr}}

\begin{abstract}
In the paper we consider $T_{1},\dots, T_{d}$ absolute
contractions of von Neumann algebra $\M$ with normal, semi-finite,
faithful trace, and prove that for every bounded Besicovitch
weight $\{a(\kb)\}_{\kb\in\bn^d}$ and every $x\in L_{p}(\M)$,
($p>1$) the averages
\begin{equation*}
A_{\Nb}(x)=\frac{1}{|\Nb|}\sum\limits_{\kb=1}^{\Nb}a(\kb)\Tb^{\kb}(x).
\end{equation*} converge bilaterally almost uniformly in
$L_{p}(\M)$.
\end{abstract}
\maketitle

\footnotetext[1]{ Current address (F.M.): Department of Comput. \&
Theor. Sci., Faculty of Sciences, IIUM, P.O. Box, 141, 25710,
Kuantan, Pahang, Malaysia }

\section{Introduction}

It is known the almost everywhere convergence of sequences of
operators were applied to study of the individual ergodic theorem
in von Neumann algebras by many authors
\cite{g},\cite{gg},\cite{La},\cite{ye} (see \cite{ja1,ja2} for
review). Very recently, in \cite{jx1,jx2} various maximal ergodic
theorems in noncommutative $L_p$-spaces have been established. As
an application of such results  the corresponding individual
ergodic theorems were obtained. Study of the almost everywhere
convergence of weighted averages in von Neumann algebras is
relatively new. In \cite{he}, the Besicovitch weighted ergodic
theorem were firstly proved in a semifinite von Neumann algebra
with a faithful normal state, which was a generalization of
\cite{ry} to a noncommutative setting. In \cite{gl} the
non-commutative Banach principle firstly obtained, further using
it in \cite{lm} a particular case of one-dimensional Besicovitch
weighted ergodic theorems were proved in the space of integrable
operators affiliated with a von Neumann algebra. Latter on, in
\cite{cls} by means of the that Principle the Besicovitch weighted
ergodic theorem has been proved in noncommutative $L_1$-spaces.

In \cite{jo} the multiparameter Besicovitch weights were
introduced and weighted ergodic theorems were obtained in
commutative $L_p$-spaces. Some other related investigations were
done in \cite{b, djo, lo}. The present paper is devoted to the
noncommutative extension of that result. Further,  we are going to
prove the bilateral almost uniform convergence of weighted
multiparameter averages with respect to bounded Besicovitch
families in noncommutative $L_{p}$-spaces. To prove it, we use the
maximal ergodic inequality for absolute contractions given in
\cite{jx2}.

\section{Preliminaries and Notations}

In what follows, $\M$ would be a semifinite von Neumann algebra
equipped with a normal semifinite faithful trace $\t$. Let $S_+$
denote the set of all $x\in \M_+$ such that $\t({\rm
supp}\,x)<\infty$, where ${\rm supp}\,x$ denotes the support of
$x$. Let $S$  be the linear span of $S_+$. Then $S$ is a
$w*$-dense $\ast$-subalgebra of $\M$. Given $1\leq p<\infty$, we
define
 $$\|x\|_p=\big[\t(|x|^p)\big]^{1/p},\quad x\in S,$$
where $|x|=(x^*x)^{1/2}$ is the modulus of $x$. Then
$(S,\;\|\cdot\|_p)$ is a normed space, whose completion is the
noncommutative $L_p$-space associated with $(\M, \t)$, denoted by
$L_p(\M,\t)$ or simply by $L_p(\M)$. As usual, we set
$L_\infty(\M,\t)=\M$ equipped with the operator norm
$\|\cdot\|_{\infty}$. We refer a reader to \cite{px} for more
information about noncommutative integration and to \cite{t} for
general terminology of von Neumman algebras.

Now recall some notions about the noncommutative
$L_p(\M;\ell_\infty(\bn^d))$-spaces. Let $d\geq 1$. Given $1\le
p\le\infty$, the space $L_p(\M;\ell_\infty(\bn^d))$ is defined as
the space of all families $x=(x_\kb)_{\kb\in \bn^d}$ in $L_p(\M)$
which admit a factorization of the following form: there are $a,
b\in L_{2p}(\M)$ and $y=(y_\kb)\subset L_\infty(\M)$ such that
 $$x_\kb=ay_\kb b,\quad\forall\; \kb\in \bn^d.$$
We then define
 $$\|x\|_{L_p(\M;\ell_\infty(\bn^d))}=\inf\big\{\|a\|_{2p}\,
 \sup_{\kb\in \bn^d}\|y_\kb\|_\infty\,\|b\|_{2p}\big\} ,$$
where the infimum runs over all factorizations as above. Then
$\big(L_p(\M;\ell_\infty(\bn^d)),\;
\|\cdot\|_{L_p(\M;\ell_\infty(\bn^d))}\big)$ is a Banach space
\cite{j},\cite{dj}. Note that these spaces were firstly introduced
in \cite{pis}, when $\M$ was a hyperfinite von Neumman algebra. In
\cite{j} it was shown that a family of positive elements
$x=(x_\kb)_{\kb\in\bn^d}$ belongs to $L_p(\M;\ell_\infty(\bn^d))$
iff there is $a\in L_p(\M)_+$ such that $x_\kb\le a$ for all
$\kb\in\bn^d$, and moreover,
 $$\|x\|_{L_p(\M;\ell_\infty(\bn^d))}=\inf\big\{\|a\|_{p}\;:\; a\in
 L_p(\M)_+
 \;\mbox{s.t.}\; x_\kb\le a,\ \forall\;\kb\in\bn^d\}.$$

The norm of $x$ in $ L_p(\M;\ell_\infty(\bn^d))$ will be very
often denoted by $\big\|\sup_\nb^+ x_\nb\big\|_p\ .$

We should be noted that $\big\|\sup_\nb^+ x_\nb\big\|_p$ is just a
notation for $\sup_\nb x_\nb$ does not make any sense in the
noncommutative setting. Some elementary properties of these spaces
were presented in \cite{dj,j,jx2}.

For $\nb=(n_1,\dots,n_d)\in\bn^d$ denote
$m(\nb)=\min\{n_1,\dots,n_d\}$, $M(\nb)=\max\{n_1,\dots,n_d\}$.
Let us denote $\L_{[m,n]}=\{\kb=(k_1,\cdots k_d)\in\bn^d: m\leq
m(\kb), M(\kb)\leq n\}$. In the sequel we will deal with the
following convergence, namely, a family  $(x_\nb)_{\nb\in\bn^d}$
in a Banach space $X$ converges to $x\in X$ if $\forall\e>0$,
there exists $N_0\in\bn$ such that $\|x_\nb-x\|_X<\e$ for all
$\nb:~ m(\nb)\geq N_0$.

We denote by $L_p(\M;\ell_\infty^{\L_{[m,n]}}(\bn^d))$ the
subspace of $L_p(\M;\ell_\infty(\bn^d))$ consisting of all finite
sequences $\{x_\kb\}_{\kb\in\bn^d}$ such that $ x_\kb=0$ if
$\kb\notin\L_{[m,n]}$. In accordance with our preceding
convention, the norm of $x$ in
$L_p(\M;\ell_\infty^{\L_{[m,n]}}(\bn^d))$ will be denoted by
$\|\sup_{\kb\in\L_{[m,n]}}^+x_\kb\|_p$.

Now introduce a subspace $L_p(\M; c_0(\bn^d))$ of $L_p(\M;
\ell_\infty(\bn^d))$, which is defined as the space of all
families $(x_\nb)_{\nb\in\bn^d}\subset L_p(\M)$ such that there
are $a,b\in L_{2p}(\M)$ and $(y_\nb)\subset\M$ verifying
 $$x_\nb=ay_\nb b\quad\mbox{and}\quad \lim_{\nb\to\infty} \|y_\nb\|_\infty=0.$$
One can check that $L_p(\M; c_0(\bn^d))$ is a closed subspace of
$L_p(\M; \ell_\infty(\bn^d))$ (see \cite{dj,j}, for more details)
and
 $$\big\|{\sup_\nb}^+ x_\nb\big\|_p=\inf\big\{\|a\|_{2p}\,
 \sup_{\nb\in\bn^d}\|y_\nb\|_\infty\,\|b\|_{2p}\big\},$$
where the infimum runs over all factorizations of $(x_\nb)$ as
above.

Let $\M$ be, as before, a von Neumann algebra equipped with a
semifinite normal faithful trace $\t$. Let $x,(x_\nb)\subset
L_p(\M)$. A family $(x_\nb)$ is said to converge {\it bilaterally
almost uniformly} $(${\it b.a.u.} in short$)$ to $x$ if for every
$\e>0$ there is a projection $e\in \M$ such that
 $$\t(e^\perp)<\e\quad \mbox{and}\quad
 \lim_{m(\nb)\to\infty}\|e(x_\nb-x)e\|_\infty=0.$$

In the commutative case, the convergence in the definition above
is equivalent to the usual almost everywhere convergence by virtue
of Egorov's theorem \cite{pa,si}.

The following Lemma gives a relation between the b.a.u.
convergence and $L_p(\M; c_0(\bn^d))$.

\begin{lem}\label{au}
If for every $\pb\in\bn^d$ $\{x_{\nb+\pb}-x_\nb\}_{\nb\in\bn^d}\in
L_p(\M; c_0(\bn^d))$ with $1\le p<\infty$, then $x_\nb$ converges
b.a.u. to some $x$ from $L_p(\M)$.
\end{lem}

The proof immediately follows from Lemma 6.2 \cite{jx2}, Theorems
1.2 and 2.3 \cite{cls}. \\

Let $(Z,\F,\mu)$ be a measurable space with a probability measure
$\mu$. Let $\widetilde \M$ be the von Neumann algebra of all
essentially bounded ultra-weakly measurable functions
$h:(Z,\mu)\to \M$ equipped with the trace
$$\tilde \tau(h)=\int_{Z}\tau(h(z))d\mu(z),$$
and let $\widetilde {L_{p}}=L_{p}(\widetilde\M,\tilde \tau)$. It
is known \cite{bgm} that the space $\widetilde {L_{p}}$ is
isomorphic to $L_{p}(Z,\mu;L_p(\M))$.

For the sake of completeness, we provide the proof for the next
lemma, which is an analog of Lemma 2 in \cite{da} (see also
\cite{cls}).

\begin{lem}\label{dec}  Let $\{x_{\nb}\}_{\nb\in\bn^d}\in L_{p}(\widetilde{\M},
c_{0}(\bn^d))$. Then $\{x_{\nb}(z)\}_{\nb\in\bn^d}\in L_{p}(\M,
c_{0}(\bn^d))$ for almost all $z\in Z$.
\end{lem}

\begin{proof} From the definition of $L_{p}(\widetilde{\M},
c_{0}(\bn^d))$ we have $x_{\nb}=ay_{\nb}b$, where $a, b\in
L_{2p}(Z,\mu; L^{p}(\M))$, $y_{\nb}\in\widetilde\M$, and
$\|y_{\nb}\|_{\widetilde\M}\rightarrow 0$ as
$\nb\rightarrow\infty.$

For any $z\in Z$ consider $x_{\nb}(z)=a(z)y_{\nb}(z)b(z)$. We have
$$\|y_{\nb}\|_{\widetilde{M}}=ess\sup_{z\in Z}\|y_{\nb}(z)\|_{\M}\rightarrow
0,$$ therefore $\|y_{\nb}(z)\|_{\M}\rightarrow 0$, $\forall z\in
Z\backslash D,$ where $\mu(D)=0.$

Since $a, b\in L_{2p}(Z, L_{p}(\M))$ we conclude that

$$\underset{Z}{\int}\|a(z)\|_{L_{p}(\M)}^{2p}d\mu<\infty;
\quad \underset{Z}{\int}\|b(z)\|_{L_{p}(\M)}^{2p}d\mu<\infty$$

Denote
$$
Z_{k}^{(\phi)}\small=\{z\in Z
:\|\phi(z)\|_{L_{p}(\M)}^{2p}>k^{2}\},$$ and put
$Z^{(\phi)}=\bigcap\limits_{n=1}^{\infty}\bigcup\limits_{k\geq
n}\large Z_{k}^{(\phi)}$, where $\phi=a, b$.  Then
\begin{eqnarray*}
\mu(\large Z^{(\phi)})&\leq &\mu\bigg(\bigcup\limits_{k\geq n}
Z_{k}^{(\phi)}\bigg)\\
&\leq&\sum\limits_{k=n}^{\infty}\mu(\large
Z_{k}^{(\phi)})\\
&\leq&\|\phi\|^{1/2p}_{L_{2p}(Z,
L_{p}(\M))}\sum\limits_{k=n}^{\infty}\frac{1}{k^{2}}
\underset{n\rightarrow\infty}\longrightarrow 0.
\end{eqnarray*}

Here, we have used that
$$\mu(\large
Z_{k}^{(\phi)})\leq\frac{1}{k^{2}}\int\|\phi(z)\|_{L_{p}(\M)}^{2p}d\mu$$

So, $\mu(\large Z^{(\phi)})=0$. Putting $N=\large Z^{(a)}\bigcup
Z^{(b)}\bigcup D$, we have $\m(N)=0$, and $a(z),b(z)\in
L_{2p}(\M)$, $y_{\nb}(z)\in\M$ with $\|y_\nb(z)\|_\M\to 0$ as
$\nb\to\infty$ for every $z\in Z\setminus N$.  Hence
$\{x_{\nb}(z)\}_{\nb\in\bn^d}\in L_{p}(\M, c_{0}(\bn^d))$ a.e.
$z\in Z$.
\end{proof}

Let $T: \M\to\M$ be a linear map. We say that $T$ is an {\it
absolute contraction} if it satisfies the following conditions
 \begin{itemize}
 \item[(i)]  $T$ is a contraction on $\M$:
$\|Tx\|_\infty\le \|x\|_\infty$ for all $x\in \M$. \item[(ii)] $T$
is positive: $Tx\ge 0$ if $x\ge 0$. \item[(iii)] $\t\circ T\le\t$:
$\t(T(x))\le\t(x)$ for all $x\in L_1(\M)\cap\M_+$.
 \end{itemize}

It is well known \cite{jx2}, \cite{ye} that if $T$ satisfies these
properties, then $T$ naturally extends to a contraction on
$L_p(\M)$ for all $1\le p<\infty$.

In the sequel, unless explicitly specified otherwise, $T$ will
always denote an absolute contraction of $\M$. The same symbol $T$
will also stand for the extensions of $T$ on $L_p(\M)$. Now let
$T_1,\dots, T_d$ be such kind of mappings.

We form their ergodic averages: \begin{equation}\label{erg}
M_{\Nb}(\Tb)=\frac{1}{|\Nb|}\sum_{\kb=1}^{\Nb}\Tb^{\kb}
\end{equation}
where $\Tb^{\kb}=T_d^{k_d}\cdots T_1^{k_1}$ with
$\Nb=(N_1,\dots,N_d)$, $\kb=(k_1,\dots, k_d)$ and $|\Nb|=N_1\cdots
N_d$.

In \cite{jx2} the following maximal inequality has been proved

\begin{thm}\label{max}
Let $p>1$ and $T_1,\dots\,,T_d$ be absolute contractions of $\M$.
For any $x\in L_p(M)_+$ there is $a\in L_p(M)_+$ such that
\begin{equation*}
M_{\Nb}(\Tb)(x)\leq a, \ \ \forall \Nb\in\bn^d \ \ \textrm{and} \
\ \|a\|_p\leq C^d_p\|x\|_p,
\end{equation*}
where $C_p$ is a positive constant depending only on p. Moreover,
$C_p\le C\, p^2(p-1)^{-2}$  and this is the optimal order of $C_p$
as $p\to 1$.
 \end{thm}

In the present paper we are going to consider the Besicovich
weights in $\bn^d$. It is said \cite{jo} that a family of complex
numbers $\{a(\kb)\}$ to be {\it Besicovitch weight} if for every
$\e> 0$ there is a family of trigonometric polynomials
$\{P_\e(\kb)\}$ in $d$ variables such that
\begin{equation}\label{B}
\limsup_{m(\Nb)\to\infty}\frac{1}{|\Nb|}\sum_{\kb=1}^{\Nb} |a(\kb)
- P_\e(\kb)| <\e.
\end{equation}

A Besicovitch weight $\{a(\kb)\}$ is said to be {\it bounded} if
$\{a(\kb)\}\in\ell_\infty(\bn^d)$. Note that certain properties of
Besicovich weights were studied in \cite{jo,lo, ry}.

\section{ Weighted vector ergodic theorem for noncommutative $L_{p}-$ spaces}

In this section we are going to prove the following a weighted
ergodic theorem for absolute contractions.

\begin{thm}\label{WET} Let $T_{1},\dots, T_{d}$ be absolute contractions of von Neumman algebra $\M$. For
every bounded Besicovitch weight $\{a(\kb)\}$ and every $x\in
L_{p}(\M)$ ($p>1$) the averages
\begin{equation}\label{werg}
A_{\Nb}(x)=\frac{1}{|\Nb|}\sum\limits_{\kb=1}^{\Nb}a(\kb)\Tb^{\kb}(x).
\end{equation}
converge b.a.u. in $L_{p}(\M)$.
\end{thm}

To prove the main result now we need some auxiliary Lemmas. Let us
first recall the following well-know principle:

\begin{thm} Let $X,Y$ be Banach spaces and $Z$ be a subspace of
$Y$. Assume that $T:X\to Y$ is a linear continuous  mapping such
that $T(X_0)\subset Z$ for a dense subset $X_0\subset X$. Then
$T(X)\subset Z$.
\end{thm}

From this principle, by taking $Z= L_{p}(\M, c_{0}(\bn^d))$, $Y=
L_{p}(\M, \ell_{\infty}(\bn^d))$ we immediately get the following

\begin{lem}\label{BP} Let $X$ be a Banach space and $a_{\nb}:X\rightarrow
L_{p}(\M, \ell_{\infty}(\bn^d))$ $(\nb\in\bn^d)$ be linear
mappings such that
\begin{enumerate}
\item[(i)] $\|a_{\nb}(x)\|_{L_{p}(\M, \ell_{\infty}(\bn^d))}\leq
C\|x\|_{X},\quad \forall x \in X$; \item[(ii)]
$\{a_{\nb+\pb}(x)-a_{\nb}(x)\}_{\nb\in\bn^d}\in L_{p}(\M,
c_{0}(\bn^d))$ for every $\pb\in\bn^d$ and $x\in X_{0}$, where
${X_{0}}$ is a dense subspace of $X$.
\end{enumerate}
Then $\{a_{\nb+\pb}(x)-a_{\nb}(x)\}_{\nb\in\bn^d}\in L_{p}(\M;
c_{0}(\bn^d))$ for all $x\in X$.
\end{lem}

\begin{lem}\label{3} Let $T_{1}, ..., T_{d}$ be as in Theorem \ref{WET}. Then for every trigonometric
polynomial $P(\kb)$ on $\bn^d$ and every $x\in L_{p}(\M)$ the
averages \begin{equation}\label{tri}
\widetilde{A}_{\Nb}(x)=\frac{1}{|\Nb|}\sum\limits_{\kb=1}^{\Nb}P(\kb)\Tb^{\kb}(x)
\end{equation}
satisfy the following relation $\{\widetilde{A}_{\nb+\pb}(x)-
\widetilde{A}_{\nb}(x)\}_{\nb\in\bn^d}\in L_{p}(M, c_{0}(\bn^d))$
for every $\pb\in\bn^d$.
\end{lem}

\begin{proof} Let  $\mathbb{B}=\{z\in\mathbb{C}: |z|=1\}$ be the unite circle in $\bc$ with
the normalized Lebesgue measure $\s$. Let $d$ be a fixed positive
integer, then denote
\begin{equation*}
\mathbb{K}=\underbrace{\bb\times\ldots\times\bb}_{d}, \ \
\mu=\underbrace{\s\otimes\ldots\otimes\s}_d.
\end{equation*}
Now consider $\widetilde{L}_{p}=L_{p}(\widetilde\M)$, where
$\widetilde{\M}=\M\otimes L_{\infty}(\mathbb{K}, \mu)$ and
$\widetilde{\tau}=\tau\otimes\mu.$

For any $\gb=(g_{i})_{i=1}^{d},\zb=(z_{i})_{i=1}^{d}\in\mathbb{K}$
we put
$$
\gb\circ\zb=(g_1 z_{1},\dots, g_dz_{d}).
$$
Now fix $\gb\in\mathbb{K}$ and define a linear mapping
$\widetilde{\Tb}_{\gb}:\widetilde{L}_{p}\rightarrow\widetilde{L}_{p}$
by
\begin{equation}\label{abs}
\widetilde{\Tb}_{\gb}(f)(\zb)=\Tb(f(\gb\circ\zb)),
\end{equation}
where $f=f(\zb)\in \widetilde L_p$, $\zb\in\mathbb{K}$. One can
see that the mapping $\widetilde T_\gb$ is an absolute
contraction.  Therefore, according to Theorem 6.6 \cite{jx2} for
the averages $M_{\Nb}(\widetilde\Tb_{\gb})$ we have
$\{M_{\nb+\pb}(\widetilde\Tb_{\gb})(f)-M_{\nb}(\widetilde\Tb_{\gb})(f)\}_{\nb\in\bn^d}\in
L_{p}(\widetilde{\M}, c_{0}(\bn^d))$ for every $\pb\in\bn^d$ and
$f\in\widetilde L_p$. Hence, Lemma \ref{dec} implies that
$\{M_{\nb+\pb}(\widetilde\Tb_\gb)(f(\zb))-M_{\nb}(\widetilde\Tb_\gb)(f(\zb))\}_{\nb\in\bn^d}\in
L_{p}({\M}, c_{0}(\bn^d))$ a.e. $\zb\in\bk$. Now applying the
letter one to the function $f_x(\zb)= \Pi(\zb)x$, where
$\Pi(\zb)=z_{1}\cdots z_{d}$, $x\in L_p(\M)$ is fixed, which
clearly belongs to $\widetilde L_p$,  we obtain
$$
\bigg\{\Pi(\zb)\bigg(\frac{1}{|\nb+\pb|}\sum\limits_{\kb=1}^{\nb+\pb}
\Pi(\gb^{\kb})\Tb^{\kb}(x)-\frac{1}{|\nb|}\sum\limits_{\kb=1}^{\nb}
\Pi(\gb^{\kb})\Tb^{\kb}(x)\bigg)\bigg\}_{\nb\in\bn^d}\in L_{p}(\M,
c_{0}(\bn^d))
$$
for almost all $\zb\in\bk$. Here we have used \eqref{abs} to get
$(\widetilde{\Tb}_{\gb})^\kb(f_x)(\zb)=\Pi(\gb^\kb)\Pi(\zb)
\Tb^\kb x$ for every $\kb\in\bn^d$, where
$\gb^{\kb}=(g_{1}^{k_1},\dots, g_{d}^{k_d})$.

Consequently, due to $\Pi(\zb)\neq 0$, we conclude that
$$\bigg\{\frac{1}{|\nb+\pb|}\sum\limits_{\kb=1}^{\nb+\pb}
\Pi(\gb^{\kb})\Tb^{\kb}(x)-\frac{1}{|\nb|}\sum\limits_{\kb=1}^{\nb}
\Pi(\gb^{\kb})\Tb^{\kb}(x)\bigg\}_{\nb\in\bn^d}\in L_{p}(\M,
c_{0}(\bn^d)).$$

One can see that the theorem holds for finite linear combinations
of $\Pi(\gb^{\kb})$, hence holds for trigonometric polynomials in
$d$ variables.
\end{proof}

Now let us turn to the proof of the main Theorem \ref{WET}.

\begin{proof} Take $x\in L_{p}(\M)_+$. Without lost of generality we may assume that
$|a(\kb)|\leq 1$ for all $\kb\in\bn^d$. We are going to show that
$A_{\Nb}(x)$ belongs to $L_{p}(\M, \ell_{\infty}(\bn^d))$.
 To do it let us consider \footnote{Note
that $\{a(\kb)\}$ are complex numbers, therefore we need to
consider their real and imaginary parts separately.}
\begin{eqnarray}\label{AR}
&& A_{\Nb}^{(R)}(x)=\frac{1}{|\Nb|}\sum\limits_{\kb=1}^{\Nb}
\Re(a(\kb))\Tb^{\kb}(x),\\\label{AI} &&
A_{N}^{(I)}(x)=\frac{1}{|\Nb|}\sum\limits_{\kb=1}^{\Nb}
\Im(a(\kb))\Tb^{\kb}(x).
\end{eqnarray}

Due to our assumption (i.e. $|a(\kb)|\leq 1$) and Theorem
\ref{max} one has
\begin{eqnarray}\label{AR1}
-a\leq-M_{\Nb}(\Tb)(x)\leq A_{\Nb}^{(R)}(x)\leq
M_{\Nb}(\Tb)(x)\leq a,
\end{eqnarray}
where $\|a\|_p\leq C_p^d\|x\|_p$. From \eqref{AR1} we have
$A_{\Nb}^{(R)}(x)\in L_{p}(\M, \ell_{\infty}(\bn^d))$, moreover
\begin{eqnarray}\label{sup}
&&\|{\sup_\Nb}^+ A_{\Nb}^{(R)}(x)\|_{L_{p}(\M,
\ell_{\infty}(\bn^d))}\leq C\|x\|_{p},
 \end{eqnarray}
for some constant $C$ (more exactly, $C=3C_p^d$). Similarly, one
gets
\begin{eqnarray}\label{sup1}
 &&\|{\sup_\Nb}^+A_{\Nb}^{(I)}(x)\|_{L_{p}(\M,
\ell_{\infty}(\bn^d))}\leq C\|x\|_{p}.
\end{eqnarray}

Consequently, \eqref{sup},\eqref{sup1} imply
\begin{equation}\label{max4}
\|{\sup_\Nb}^+A_{\Nb}(x)\|_{L_{p}(\M, \ell_{\infty}(\bn^d))}\leq
2C\|x\|_{p},\quad \forall x\in L_{p}(\M)_+.
\end{equation}

Take any $x\in L_p(\M)$, then one can be represented as
$x=\sum_{k=0}^3i^kx_k$, where $x_k\in L_p(\M)_+$, $k=0,1,2,3$.
Therefore, the inequality \eqref{max4} implies that
\begin{equation}\label{max41}
\|{\sup_\Nb}^+ A_{\Nb}(x)\|_{L_{p}(\M, \ell_{\infty}(\bn^d))}\leq
8C\|x\|_{p},\quad \forall x\in L_{p}(\M).
\end{equation}

Now let us assume that $x\in L_1(\M)\cap\M$, then $x\in L_q(\M)$
for any $1<q<\infty$. For an arbitrary $\e>0$ due to the
definition of the Besicovitch weight there is a trigonometric
polynomial $\{P_\e(\kb)\}$ for which  \eqref{B} holds. Let
$\widetilde{A}_{\Nb}(x)$ be the corresponding averages (see
\eqref{tri}). Then from \eqref{max41} we conclude that
$\{A_{\Nb}(x)\}\in L_{q}(\M, \ell_{\infty}(\bn^d))$,
$\{\widetilde{A}_{\Nb}(x)\}\in L_{q}(\M, \ell_{\infty}(\bn^d))$,
so $\{A_{\Nb }(x)-\widetilde{A}_{\Nb}(x)\}\in L_{q}(\M,
\ell_{\infty}(\bn^d))$. Moreover,
\begin{eqnarray*}\label{max50}
\|A_{\kb}(x)-\widetilde{A}_{\kb}(x))\|_{\infty}&\leq&
\frac{1}{|\Nb|}\sum_{\kb=1}^{\Nb}|a(\kb)-P_\e(\kb)|\|\Tb^\kb
(x)\|_\infty \nonumber\\
&\leq&\|x\|_\infty\frac{1}{|\Nb|}\sum_{\kb=1}^\Nb|a(\kb)-P_\e(\kb)|<\e\|x\|_\infty
\end{eqnarray*}

From the last relation and Proposition 2.5 \cite{jx2} one finds
\begin{eqnarray}\label{max5}
\|{\sup_{\kb\in\L_{[k,n]}}}^+(A_{\kb}(x)-
\widetilde{A}_{\kb}(x))\|_{L_{p}(\M, \ell_{\infty}(\bn^d))}&\leq&
\sup_{\kb\in\L_{[k,n]}}\|A_{\kb}(x)-\widetilde{A}_{\kb}(x))\|_{\infty}^{1-\frac{q}{p}}\nonumber\\
&&\times\|{\sup_{\kb\in\L_{[k,n]}}}^+(A_{\kb}(x)-\widetilde{A}_{\kb}(x))\|_{L_{q}(\M, \ell_{\infty}(\bn^d))}^{\frac{q}{p}}\nonumber\\
&\leq&
\e^{1-\frac{q}{p}}\|x\|_{\infty}^{1-\frac{q}{p}}\|{\sup_{\kb\in\bn^d}}^+(A_{\kb}(x)-\widetilde{A}_{\kb}(x))
\|_{L_{q}(\M, \ell_{\infty}(\bn^d))}^{\frac{q}{p}}
\end{eqnarray}

Define a sequence $b^{(k)}=(b^{(k)}_{\kb})_{\kb\in\bn^d}\in
L_p(\M,c_0(\bn^d))$ as follows:
\begin{equation*}
b^{(k)}_\kb= \left\{ \begin{array}{ll}
A_{\kb}(x)-\widetilde{A}_{\kb}(x),  \ \ \textrm{if} \ \ \kb\in\L_{[0,k]},\\
0, \qquad \qquad \qquad \ \  \ \textrm{if} \ \
\kb\notin\L_{[0,k]}.
\end{array}
\right.
\end{equation*}

From \eqref{max5} one gets that
$b^{(k)}\rightarrow\{A_{\kb}(x)-\widetilde{A}_{\kb}(x)\}$ in
$L_{p}(\M, \ell_{\infty}(\bn^d))$ as $k\to\infty$. Since
$\{b^{(k)}\}\subset L_{p}(\M, c_{0}(\bn^d))$ we obtain $\{
A_{\kb}(x)-\widetilde{A}_{\kb}(x)\}\in L_{p}(\M, c_{0}(\bn^d))$.

According to Lemma \ref{3} we have already known  that
$\{\widetilde{A}_{\nb+\pb}(x)-\widetilde{A}_{\nb}(x)\}_{\nb\in\bn^d}\in
L_{p}(\M, c_{0}(\bn^d))$. Consequently, the equality
$$
A_{\nb+\pb}(x)-{A}_{\nb}(x)=(A_{\nb+\pb}(x)-\widetilde{A}_{\nb+\pb}(x))+(\widetilde{A}_{\nb}(x)-A_{\nb}(x))+
(\widetilde{A}_{\nb+\pb}(x)-\widetilde{A}_{\nb}(x))
$$
implies that $\{A_{\nb+\pb}(x)-{A}_{\nb}(x)\}_{\nb\in\bn^d}\in
L_{p}(\M, c_{0}(\bn^d)).$

Now by means of the density of  $L_1(\M)\cap\M$ in $L_p(\M)$ and
\eqref{max41} with Lemma \ref{BP}, we have
$\{A_{\nb+\pb}(x)-{A}_{\nb}(x)\}_{\nb\in\bn^d}\in L_{p}(\M,
c_{0}(\bn^d))$ for any $x\in L_p(\M)$. Hence, Lemma \ref{au}
implies that the required assertion.
\end{proof}

{\bf Remark.} 1. The proved theorem extends the results of the
papers \cite{cls} to multiparameter. When $a(\kb)\equiv 1$, then
we get an extension of  \cite{sk},\cite{pe} to $L_p$-spaces.

2.  Similar results for multiparameter Besicovitch weights in a
commutative setting were obtained  in \cite{b,jo}.

\section*{acknowledgments} The first named author (F.M.) is partially supported by
FCT grant SFRH/BPD/ 17419 /2004. The authors also would like to
thank to the referee for his useful suggestions which allowed us
to improve the text of the paper.

\end{document}